\documentclass[12pt,letterpaper]{amsart}
\usepackage{amssymb}
\usepackage{amsmath}
\usepackage{amsthm}
\usepackage{overpic}
\usepackage{epstopdf}
\usepackage{fullpage}
\usepackage{hyperref}

\begin{document}

\theoremstyle{plain}
\numberwithin{equation}{section}
\newtheorem{theorem}[equation]{Theorem}
\newtheorem{proposition}[equation]{Proposition}
\newtheorem{lemma}[equation]{Lemma}
\newtheorem{corollary}[equation]{Corollary}
\theoremstyle{definition}
\newtheorem{definition}[equation]{Definition}
\newtheorem{remark}[equation]{Remark}
\renewcommand{\thefigure}{\theequation}

\newcommand{\R}{\mathbb{R}}
\newcommand{\Z}{\mathbb{Z}}
\newcommand{\N}{\mathbb{N}}
\newcommand{\I}{\mathcal{I}}
\renewcommand{\O}{\mathcal{O}}
\renewcommand{\S}{\mathcal{S}}
\newcommand{\tildeSigma}{\tilde{\Sigma}}
\newcommand{\iso}{\text{iso}}
\newcommand{\im}{\text{im}}
\newcommand{\<}{\left\langle}
\renewcommand{\>}{\right\rangle}

\title{Examples of Stable Embedded Minimal Spheres without Area Bounds}
\author{Joel Kramer}

\maketitle


\section{Introduction}

	Throughout this paper, we use the $C^2$-topology on the space of metrics on a manifold. The main result will be the following theorem:
		\begin{theorem}\label{maintheorem}
	Let $M^3$ be a three-manifold. There exists an open, nonempty set of metrics on $M$ for each of which there are stable embedded minimal two-spheres of arbitrarily large area. 
		\end{theorem}
	This answers affirmatively a question posed in \cite{coldingminicozzi2004}. In \cite{coldingminicozzi1999b}, T. H. Colding and W. P. Minicozzi II proved an analogous theorem, stating that there exists an open set of metrics on $M$ such that there are embedded minimal tori of arbitrarily large area. This result was extended to surfaces of positive genus in \cite{dean2003} by B. Dean. The genus zero case has remained an open problem since then, largely due to the fact that the fundamental group of a sphere is trivial. 
	
	The approaches to building the positive genus surfaces followed a simple structure. One considers a closed three-dimensional submanifold of the unit three-ball which is thought of as a subset of the larger manifold. A sequence of surfaces is then constructed such that the area, allowing for certain variance of the surface, seems to approach infinity. Given any element of this sequence, a result by R. Schoen and S. T. Yau (see \cite{schoenyau1979}) is used to show that there is a stable minimal surface which closely approximates this element. The final step is to show that the area becomes unbounded.

	The result of Schoen and Yau mentioned above states that given a closed non-simply connected embedded surface with  for which the inclusion map induces an injection of the fundamental group, one can always find a closed embedded stable minimal surface of the same genus whose fundamental group has the same image under the induced map from inclusion. Replacing this is a result by W. Meeks III, L. Simon, and S. T. Yau (see \cite{meekssimonyau1982}) in which one can find a set of closed stable embedded minimal surfaces which is obtained from the original by pinching off parts of the surface and varying each isotopically. This is discussed in detail in section (\ref{treductionsection}).

It should also be noted that J. Hass, P. Norbury, and J. H. Rubinstein constructed embedded minimal spheres of unbounded morse index in \cite{hnr2003} by methods which are significantly different than those presented here. Also, in contrast of the main theorem, H. Choi and A. Wang proved in \cite{choiwang1983} that there is an open set of metrics, namely those with positive ricci curvature, for which there is a uniform area bound on compact embedded minimal surfaces depending only on the genus and the ambient metric. The author also suspects that it may be possible to provide a bound (depending on the metric) for the genus of any embedded minimal surface in the unit ball.

	The author would like to thank Professor Minicozzi for bringing this problem to his attention and his continued guidance.


\section{Notation}

	In this paper, $I=[0,1]\subset \R$ will denote the unit interval. For a subset $U$ of a topological space $X$, $U^\circ$ will be the interior and $\overline{U}$ will be the closure of $U$ in $X$. For a continuous function $f:X\rightarrow Y$, $f_*$ will denote the map induced on homotopy and $im(f)=f(X)$ will be the image of $f$.

 	Let $N$ be an arbitrary three dimensional differentiable manifold.we will say  $U\subset N$ is \emph{homotopic} to $U'\subset N$ if there is a continuous homotopy $\varphi:I\times U\rightarrow N$ such that $\varphi_0(U)=U$ and $\varphi_1(U)=V$. A subset $U\subset N$ is said to be \emph{homotopically trivial} if $U$ is homotopic to a point.

	An \emph{isotopy} (also known as an \emph{ambient isotopy}) from a set $U\subseteq N$ to the set $U'\subseteq N$ is a continuous one parameter family of diffeomorphisms $\varphi:I\times N\rightarrow N$ such that $\varphi_0(U)=U$ and $\varphi_1(U)=U'$. $\I(U)$ will denote the class of sets in $N$ isotopic to $U$. 

	We will use the following notation to define surfaces with multiplicity.	Let $m_i\in \Z$ and $\Sigma_i\subset N$ closed. We define the sum $m_1\Sigma_1+\ldots m_R\Sigma_R$ be the set $\{(m_i,\Sigma_i):i=1\ldots R\}$. Given $\Sigma_0+m\Sigma$, if $\Sigma$ has a smooth choice of unit normal,  $\I(\Sigma_0+ (m\Sigma))$ will denote the isotopy class of the union of $\Sigma_0$ and the graphs of the constant functions $k\epsilon/m$ over $\Sigma$ for $k=1,\ldots,m$ where $\epsilon$ is arbitrarily small. We define $\mathcal{I}(\Sigma_0\cup 0\Sigma)$ to be $\mathcal{I}(\Sigma_0)$.

	If $N$ has a riemannian metric, the area of a subsurface $\Sigma$ will be denoted $|\Sigma|$.

\section{Construction of $M_{k/2}$}\label{constmn}

\begin{figure}
			\begin{center}
				\begin{overpic}[scale=1]{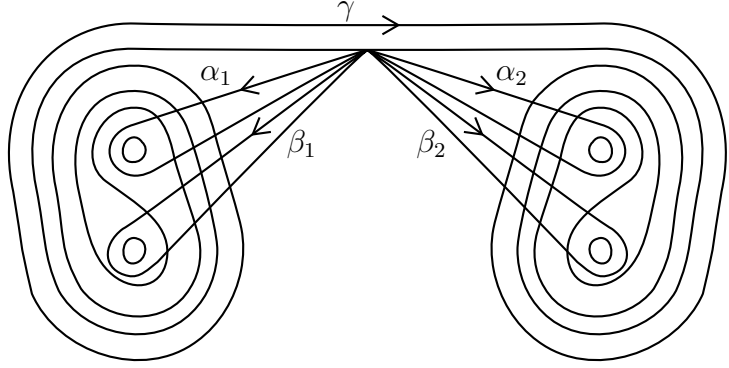}
					\put(27,39){$\alpha_1$}
					\put(68,39){$\alpha_2$}
					\put(39,29){$\beta_1$}
					\put(57,29){$\beta_2$}
					\put(46,48){$\gamma$}
				\end{overpic}
				\caption{$\alpha_1$, $\beta_1$, $\alpha_2$, $\beta_2$, and $\gamma_2$}\label{curve}
			\end{center}
		\end{figure}

	Consider the open unit ball $B_2\subset\R^2$ and let $r(x,y)=(-x,y)$ where $x,y$ are the usual Euclidean coordinates. Choose two disjoint closed sub-disks $D_1,D_2\subset\{(x,y)\in B_2^\circ:x<0\}$. Let $\Omega$ be the closure of $B_2-(D_1\cup D_2\cup r(D_1\cup D_2))$. Let $L'=\{(x,y)\in B_2:x=0\}=\{(x,y)\in \Omega:x=0\}$ and $x_0\in L'$. $\pi_1(\Omega,x_0)$ is generated by the loops $\alpha_1$, $\beta_1$, $\alpha_2=r_*(\alpha_1)$, and $\beta_2=r_*(\beta_1)$ where $\alpha_1$ and $\beta_1$ go around $D_1$ and $D_2$ counterclockwise respectively. For $n\in \Z$, define $\gamma_n,\gamma_{n+1/2}\in\pi_1(\Omega)$ as
	\begin{align}
		\begin{split}\label{gamma}
			\gamma_n 		&=(\alpha_1\beta_1)^n\alpha_1^{-1}(\alpha_1\beta_1)^{-n}(\alpha_2\beta_2)^n\alpha_2(\alpha_2\beta_2)^{-n}
		\end{split}\\
		\begin{split}\label{gammastar}
			\gamma_{n+1/2} 	&=(\alpha_1\beta_1)^{n}\alpha_1(\alpha_1\beta_1)^{-n-1}(\alpha_2\beta_2)^{n+1}\alpha_2^{-1}(\alpha_2\beta_2)^{-n}
		\end{split}
	\end{align}
	(see figure (\ref{curve})). 

		\addtocounter{equation}{1}

		For $k\in\Z$ choose an embedded $C^2$ curve $C_{k/2}\in\gamma_{k/2}$ such that $r(im(C_{k/2}))=im(C_{k/2})$. Define $N$ to be the three manifold obtained by rotating $\Omega$ about the center line $L'$. That is $N=\Omega\times\R/\sim$ where $\sim$ represents the equivalence relations
	\begin{enumerate}
		\item[i.] $(x,\theta)\sim(x,\theta')$ for $x\in L'$ and $\theta,\theta'\in\R$ 
		\item[ii.] $(x,\theta)\sim(x,\theta+2k\pi)$ for $k\in\mathbb{Z}$
		\item[iii.] $(x,\theta)\sim(r(x),\theta+(2k+1)\pi)$ for $k\in\mathbb{Z}$ 
	\end{enumerate}
	We define $M_{k/2}$ to be surface obtained by rotating $C_{k/2}$ about $L'$. Then $N$ is diffeomorphic to a solid two-sphere minus two unlinked, unknotted tori, and $M_{k/2}$ are two-spheres hooking around the tori. We let $L=\{(x,\theta)/\sim:x\in L'\}$. 
	
	We also define $p:N\rightarrow\{(x,y)\in \overline B_2:x\geq 0\}$ as $p(x,y,\theta)=(x,y)$ if $x\geq 0$ and $p(-x,y)=p\circ r(x,y)$ if $x<0$. This is well defined and continuous. We call a subset $U\subset N$ \emph{rotationally symmetric} if $U=p^{-1}(V)$ for some subset $V$ of the image of $p$. We also define a \emph{euclidean $\epsilon$-neighborhood of $U\subset N$} to be the set $\{(x\in N: d(x,U)<\epsilon\}$ where $d$ is the euclidean metric induced by this construction. Note that if $U$ is rotationally symmetric, then so is any euclidean $\epsilon$-neighborhood.

	We also define two other surfaces which will come in use later. Let $\gamma_T=\alpha_1\beta_1$ and $\gamma_S=\alpha_1\beta_1\beta_2^{-1}\alpha_2^{-1}$ be elements of $\pi_1(\Omega)$. Choose an embedded $C^2$ curve $C_T\in\gamma_T$ such that $C_T$ lies on one side of $L'$ in $\Omega$ and choose an embedded $C^2$ curve $C_S\in\gamma_S$ such that $C_S$ is symmetric with respect to reflection about $L'$. Then define $T$ and $S$ by rotating $C_T$ and $C_S$ respectively about $L'$ as described above. Conceptually, $T$ is a torus containing both of the torus components of $\partial N$ and $S$ is a sphere containing both of the tori which is homotopic to the spherical component of $\partial N$. We can assume, without loss of generality, we have chosen $M_{k/2}$ and $S$ such that $M_{k/2}\cap S=\emptyset$.  

\begin{proposition}
	 If $M_{k/2}$ is homotopic to $M_{l/2}$, then $l=k$. 
\end{proposition}

This necessarily implies isotopic distinction as well. To prove this, we will need a few higher homotopy results. Let $X$ be a path connected space. For a path $\gamma:I\rightarrow X$ with $\gamma(0)=x_0$ and $\gamma(1)=x_1$, we can define a natural change of base point map (also called $\gamma$) from $\pi_n(X,x_1)$ to $\pi_n(X,x_0)$ as follows: let $y_1,\ldots,y_n$ be coordinates for $I^n$ and let $U=\{y\in I^n:1/4\leq y_i\leq 3/4\text{ for }i=1\ldots n\}$. Choose a map $g:I^n-U\rightarrow I$ with $g(\partial I^n)=0$ and $g(\partial U)=1$. Then for $[f]\in\pi_n(X,x_1)$, let $f_\gamma : I^n\rightarrow X$ be defined as:
	\begin{equation}\label{action}
		f_\gamma(y) = 
			\begin{cases}
				\gamma\circ g(y)										&\text{if }	y\in(I^n-U)\\
						f(2(y_1-1/4),\ldots,2(y_n-1/4))	&\text{if }	y\in U
			\end{cases}
	\end{equation}
	We let $\gamma([f])=[f_\gamma]$ (note that this is homotopically independent of our choice of $g$). If $x_0=x_1$, then this represents an action of $G=\pi_1(X,x_0)$ on $\pi_n(X,x_0)$. For $A\in\pi_n(X,x_0)$, we call $GA=\{\gamma A:\gamma\in G\}$ the \emph{orbit of $A$}. We let $\mathcal{O}_n(X)$ denote the set of orbits in $\pi_n(X)$. 

	In the following lemma, we will let $\mathcal S_n(X)$ be the set of homotopy classes of maps from $S^n$ into $X$.

\begin{lemma}\label{Theta}
	Let $q:I^n\rightarrow S^n$ be the quotient map taking $\partial I^n$ to a single point and for $[f]\in\pi_n(X)$ let $\bar f: S^n\rightarrow X$ be the map such that $f=\bar f \circ q$. Let $\Theta:\pi_n(X)\rightarrow \S_n(X)$ be the map $\Theta([f])=[\bar f]$. $\Theta$ is surjective and if $\Theta(A)=\Theta(B)$, then $B=\gamma A$ for some $\gamma\in G$. That is, $\Theta$ gives a natural one to one correspondence between $\O_n(X)$ and $\S_n(X)$.
\end{lemma}
\emph{Remark:} The one dimensional analogue is the well known result that there is a one to one correspondence between conjugacy classes of $\pi_1(X)$ and homotopy classes of maps from $S^1$ into $X$.

\begin{proof}
	Let $[\bar f]\in \S^n(X)$. Let $f=\bar f\circ q$ and let $\gamma$ be a path connecting the base point $x_0$ of $\pi_n(X,x_0)$ and $f(\partial I^n)$. Then $\Theta(\gamma[f])=[\bar f]$, and surjectivity is proven.

	Assume $\Theta([f])=\Theta([g])$. Then there exists a homotopy $\bar\varphi$ such that $\bar\varphi_0=\bar f$ and $\bar\varphi_1=\bar g$. Let $\gamma:I\times I\rightarrow X$ be a one parameter family of paths defined as $\gamma_t(s)=\bar\varphi_{st}(q(\partial I^n))$. Let $\varphi$ be the homotopy $\varphi_t=(\bar\varphi_t(f))_{\gamma_t}$. This gives a homotopy between $f$ and $g_{\gamma_1}$, thus we have $[f]=\gamma_1[g]$
\end{proof}

	We will call a basis $\mathfrak{B}$ for $\pi_n(X)$ a \emph{natural basis} if for every basis element $A\in\mathfrak{B}$ and $\gamma\in\pi_1(X)$, $\gamma A\in\mathfrak{B}$. For $B\in\pi_n(X)$, define $\Phi(B)$ as the number of nonzero coefficients in the expression $B=\sum_{A\in\mathfrak{B}}\alpha_A A$. Clearly $\Phi(B)=\Phi(\gamma B)$ for $\gamma\in\pi_1(X)$, so $\Phi$ is constant on orbits. Lemma (\ref{Theta}) gives a one to one correspondence between $\O_n(X)$ and $\S_n(X)$, so we can define unambiguously $\Phi(f)=\Phi(\Theta^{-1}([f]))$ for any $f:S^2\rightarrow X$. Again, $\Phi$ is constant on homotopy classes.

	\begin{proof}[Proof of proposition]
	$N$ is homotopically equivalent to the wedge sum of two two-spheres with a line connecting the north pole of one with the south pole of the other through the connecting point. Call this space $N'$ (see figure (\ref{spheres})) and let $\phi:N\rightarrow N'$ be a homotopy equivalence. Also note that $N'$ is homotopically equivalent to the wedge sum of two two-spheres and two copies of $S^1$. This means the universal covering space of $N'$, $\tilde N'$, is homotopically equivalent to the universal covering space of $S^1\wedge S^1$ (the well known ``tree'') with two two-spheres at each vertex. Since the universal cover is simply connected, $\pi_2(\tilde{N}')$ is homomorphic to the second homology $H_2(\tilde{N}')$ which is the direct sum of countably many copies of $\mathbb{Z}$, one for each sphere. Since the $n$'th homology of a space is homomorphic to $n$'th homology of any cover for $n\geq 2$, we get that $\pi_2(N')$ is also the direct sum of countably many copies of $\mathbb{Z}$. More specifically 
\begin{equation} \pi_2(N')=\bigoplus_{\alpha\in\pi_1(N')}(\Z \oplus \Z) \end{equation}
	\addtocounter{equation}{1}\begin{figure}
			\begin{center}
				\begin{overpic}[scale=.4]{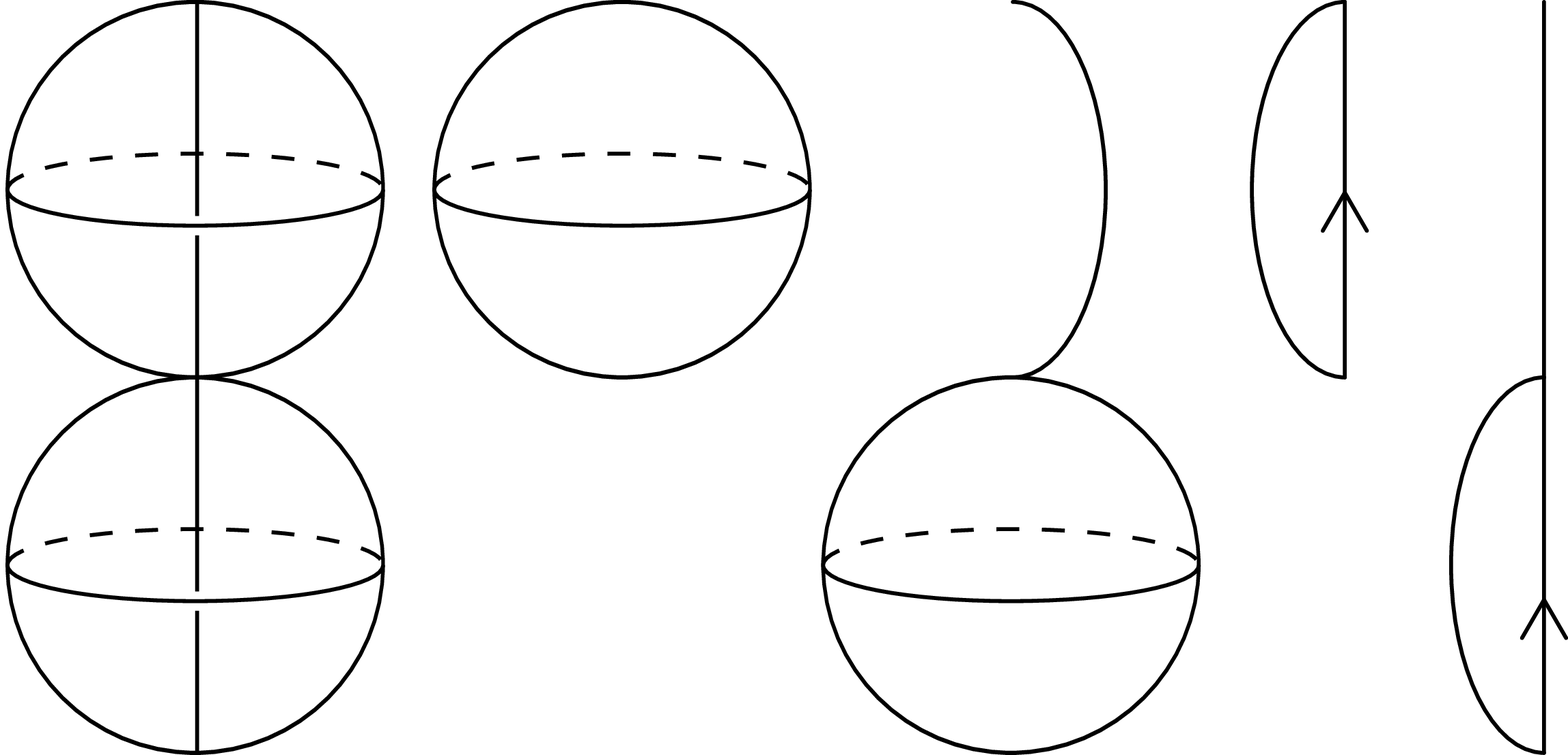}
					\put(11,49){$x_0$}
				\end{overpic}
				\caption{Respectively: $N'$, $A'_e$, $B'_e$, $\alpha'$, and $\beta'$}\label{spheres}
			\end{center}
		\end{figure}
	We can find generators of $\pi_2(N')$ as follows: Let $I^2=I\times I\subset\R^2$. Choose a base point $x_0\in N'$ as the north pole of the top sphere. Let $A'\in\pi_2(N,x_0)$ be defined by mapping $\partial I^2$ to $x_0$ and $I^2$ over the top sphere. Define $B'\in\pi_2(N,x_0)$ by mapping $\partial I^2$ to $x_0$, stringing $I^2$ around the side of the top sphere and over the bottom sphere (see figure (\ref{spheres})). For $\gamma\in\pi_1(N',x_0)$, we will define $A'_\gamma=\gamma A'$ and $B'_\gamma=\gamma B'$. We now have a set of generators for $\pi_2(N',x_0)$, namely $\left \{ A'_\gamma, B'_\gamma | \gamma\in\pi_1(N',x_0)\right \}$. For convenience, we let $C'=A'+B'$ and $C'_\gamma=\gamma C'$. We let $A_\gamma=\phi^{-1}_*(A'_\gamma)$, $B_\gamma=\phi^{-1}_*(B'_\gamma)$, and $C_\gamma=\phi^{-1}_*(C'_\gamma)$. Without loss of generality, we can assume that we have chosen the orientation of $A'$ and $B'$ such that there is an element of $C$ which is embedded in $N$.

	Let $\alpha'$ and $\beta'$ be generators of $\pi_1(N',x_0)$ such that $\alpha'$ goes around the side of the top sphere to the intersection point and back through to $x_0$ and $\beta'$ goes through the middle top sphere, around the side of the bottom sphere, and back though the center of each (see figure (\ref{spheres})). Also let $\gamma'=\alpha'\beta'$ be the path going around side of both spheres and back through the center. Again, let $\alpha=\phi^{-1}_*(\alpha')$ $\beta=\phi^{-1}_*(\beta')$ $\gamma=\phi^{-1}_*(\gamma')$.

	Let $f_{k/2}$ be the inclusion map on $M_{k/2}\subset N$. We can see that $[f_{k/2}]=\Theta(D_{k/2})$ where 
\begin{align}
	\begin{split}\label{pitwomn}
		D_n	&=C_e+C_\gamma+\ldots+C_{\gamma^{n-1}}-A_{\gamma^{n}\alpha^{-1}}-C_{\gamma^n\alpha^{-1}\gamma^{-1}}-\ldots-C_{\gamma^n\alpha^{-1}\gamma^{-n}}\\
				&=-A_{\gamma^n\alpha^{-1}}+\sum_{k=1}^n \left[C_{\gamma^{k-1}}-C_{\gamma^n\alpha^{-1}\gamma^{-k}}\right]
	\end{split}\\
	\begin{split}\label{pitwomnstar}
		D_{n+1/2}	&=C_e+C_\gamma+\ldots+C_{\gamma^{n-1}}+A_{\gamma^{n}}-C_{\gamma^n\alpha}-\ldots-C_{\gamma^n\alpha\gamma^{-n-1}}\\
							&=A_{\gamma^n}+\sum_{k=1}^{n} \left[C_{\gamma^{k-1}}-C_{\gamma^n\alpha\gamma^{-k}}\right]-C_{\gamma^n\alpha\gamma^{-n-1}}
	\end{split}
\end{align}
	for $n\in\Z$. (Note the similarity between the last term in each and equations (\ref{gamma}) and (\ref{gammastar})) So $\Phi(f_{k/2})=2k+1$ and since $\Phi$ is constant on homotopy classes, all $f_{k/2}$ are all homotopically distinct.
\end{proof}

\section{Compression and Meeks-Simon-Yau}\label{treductionsection}

	We would like to be able to find an area minimizer in the isotopy class of $M_{k/2}$. To do so we will use following result by W. Meeks III, L. Simon, and S. T. Yau. 

	Given an mean convex three-manifold $N$ and an embedded subsurface $\Sigma_0\subset N$, we can find a sequence $\Sigma_k\in\I(\Sigma_0)$ such that $|\Sigma_k|\rightarrow\inf_{\Sigma\in\I(\Sigma_0)}|\Sigma|$. The main theorem in \cite{meekssimonyau1982} shows that there is a subsequence (still called $\Sigma_k$) and compact embedded minimal surfaces $\Sigma^{(1)},\ldots,\Sigma^{(R)}$ such that
	\begin{equation}\label{surfaces}\Sigma_k\rightarrow m_1\Sigma^{(1)}+\ldots+m_R\Sigma^{(R)}\end{equation}
	in the Borel measure sense. That is, for any $f\in C_0(N)$
	\begin{equation}\lim_{k\rightarrow\infty}\int_{\Sigma_k} f =\sum_{i=1}^R m_i\int_{\Sigma^{(i)}} f\end{equation}
	Moreover, if $g_j=\text{genus}(\Sigma^{(j)})$, then we get
	\begin{equation}\label{genera}
		\sum_{j\in\mathcal{U}}\frac{1}{2}m_j(g_j-1)+\sum_{j\in\mathcal{O}}m_j g_j\leq \text{genus}(\Sigma_0)
	\end{equation}
	where $\mathcal{U}=\{j:\Sigma^{(j)} \text{ is one-sided in } N\}$ and $\mathcal{O}=\{j:\Sigma^{(j)} \text{ is two-sided in } N\}$. Also, if $\Sigma_0$ is two-sided, then each $\Sigma^{(j)}$ is stable.
	We will call $m_1\Sigma^{(1)}+\ldots+ m_R\Sigma^{(R)}$ the \emph{Meeks-Simon-Yau minimizer} or the \emph{MSY-minimizer} for $\Sigma_0$.

	To describe the relationship between $M_{k/2}$ and $\Sigma_n^{(j)}$, it is necessary to briefly cover $\gamma$-reducibility. Let $\Sigma_0$ be a compact surface without boundary (possibly with more than one connected component) embedded in a three-manifold $N^3$. For $U,V\subseteq N$, let $U\Delta V$ denote the symmetric difference $(U-V)\cup(V-U)$. 
	
	For $\gamma$ sufficiently small, we say  the surface $\tildeSigma\subseteq N$ is a $\gamma$-reduction of $\Sigma$ (denoted  $\tildeSigma\ll_\gamma \Sigma$) if 
		\begin{enumerate}
			\item $|\Sigma\Delta\tildeSigma|<2\gamma$ \label{littlearea}
			\item $\Sigma-\tildeSigma$ is diffeomorphic to an open annulus.
			\item $\tildeSigma-\Sigma$ consists of two components, each diffeomorphic to the open unit disk.
			\item\label{ball} There exists $B\subseteq N$ homeomorphic to the open unit three-ball such that $B$ is disjoint from $\Sigma\cup\tildeSigma$ and  $\partial B = \overline{\Sigma\Delta\tildeSigma}$.
			\item In the case that $\tildeSigma$ is not connected, each component is either not simply connected or else has area $>\delta^2/2$ where $\delta$ is a fixed number depending only on the injectivity radius of $N$ (see Lemma 1 of \cite{meekssimonyau1982}). \label{littlesphere}
		\end{enumerate}

	In other words, $\tildeSigma\ll_\gamma\Sigma$ if $\tildeSigma$ can be obtained from $\Sigma$ by ``pinching'' off a small part of the surface. We say $\Sigma$ is $\gamma$-irreducible if there is no $\gamma$-reduction of $\Sigma$. 
	
	\begin{proposition}\label{seq}
		Given $\Sigma_0$, $m_j$, and $\Sigma^{(j)}$ as described above, there is a pair of finite sequences, $\{\Sigma_i\}_{i=0}^{k}$ and $\{\Sigma_i'\}_{i=0}^k$, such that $\Sigma_i'\in\I(\Sigma_i)$ for all $i$, $\Sigma_{i+1}\ll_\gamma\Sigma_i'$, and $\Sigma_k'\in\I(\cup_{j=1}^R m_j\Sigma^{(j)})$. 
	\end{proposition}
	\begin{proof}
		This is a rewording of remark (3.27) of \cite{meekssimonyau1982}.
	\end{proof}
	The requirement that we have two sequences is due to (\ref{littlearea}) in the definition of $\gamma$-reduction. If we drop this restriction and modify (\ref{littlesphere}) by only requiring that the spheres be homotopically non-trivial, we get a completely topological definition as follows.  
\begin{definition}\label{tred}
	A surface $\tilde{\Sigma}$ is called a \emph{compression} of a homotopically non-trivial surface $\Sigma$, written $\tilde{\Sigma}\ll\Sigma$ if 
		\begin{enumerate}
			\item Every connected component of $\tilde{\Sigma}$ is homotopically non-constant
			\item $\Sigma-\tildeSigma$ is diffeomorphic to an open annulus.
			\item $\tildeSigma-\Sigma$ consists of two components, each diffeomorphic to the open unit disk.
			\item There exists $B\subseteq N$ diffeomorphic to the open unit three-ball such that $B$ is disjoint from $\Sigma\cup\tildeSigma$ and  $\partial B = \overline{\Sigma\Delta\tildeSigma}$.
		\end{enumerate}
	\end{definition}

We now have an analogue of proposition (\ref{seq}) for compression.

\begin{proposition}\label{Tseq}
		Given $\Sigma_0$, $m_j$, and  $\Sigma^{(j)}$ as described above, there is a finite sequence $\{\Sigma_i\}_{i=0}^{k}$ such that $\Sigma_{i+1}\ll\Sigma_i$ for $i=0,\ldots,k-1$, and $\Sigma_k\in\I(\cup_{j=1}^Rm_j\Sigma^{(j)})$.
\end{proposition}
\begin{proof}
	All that is needed to be shown here is that none of the $\gamma$-reductions in proposition (\ref{seq}) result in a surface with components which are homotopically trivial. If so, that component which is homotopically trivial is isotopic to a surface which is arbitrarily small and thus shows up as a surface with multiplicity $0$ in $\sum_{j=1}^Rm_j\Sigma^{(j)}$. Thus, this step of the reduction can be omitted without consequence.
\end{proof}

	Consider the case of $N$ being defined as in section (\ref{constmn}) with a metric defined in such a way that $N$ is mean convex, and $\Sigma_0=M_{k/2}$. Then there is a sequence $\Sigma_{k/2,l}$ for which $|\Sigma_{k/2,l}|\rightarrow \inf_{\Sigma\in\I(M_{k/2)}}|\Sigma|$ and which the sequence converges to $S_{k/2}=\sum m_{k/2,j}\Sigma_{k/2}^{(j)}$ in the Borel measure sense. Since $\Sigma_{k/2}^{(j)}$ embeds into $N$, which can be embedded into $\R^3$, we must have that each $\Sigma_{k/2}^{(j)}$ is two sided in $N$. Therefore, since $M_{k/2}$ is genus $0$, equation (\ref{genera}) tells us
	\begin{equation}\sum_{j=1}^R m_{{k/2},j} g_j \leq 0\end{equation}
	So each $\Sigma_{k/2}^{(j)}$ is an embedded stable minimal sphere. 

	For an example of compression, take $D$ to be the disk with boundary along the ``lip'' of $M_{k/2}$ and consider the surface obtained by ``cutting'' $M_{k/2}$ along $D$ (i.e., replacing the annulus $A=\{x\in M_{k/2} | d(x,\partial D)<\epsilon\}$ with $D$ moved in both normal directions by a distance of $\epsilon$). This compression is isotopic to $M_{(k-1)/2}\cup S$. As it turns out, this is one of only two isotopically distinct compressions of $M_{k/2}$. To show this, we will first show that every compression of $M_{k/2}$ can be done such that one of the components of the compression is rotationally symmetric. 

\begin{proposition}\label{onlygammared}
	If $\Sigma$ is a compression of $M_{k/2}$ then $\Sigma$ has two connected components, one of which is isotopic to $S$ and the other is isotopic to either $M_{(k+1)/2}$ or $M_{(k-1)/2}$
\end{proposition}

	To prove this, we will introduce the following notion of a surface containing another set. We will say a surface $M\subset N$ \emph{encapsulates} a set $U$ if for all $x\in U$ any path connecting $x$ to the spherical component of $\partial N$ intersects $M$ non-trivially. Note, in particular, that any surface $\Sigma_{k/2}\in\I(M_{k/2})$ encapsulates one and only one of the torus components of $\partial N$. Moreover, for all integers $n$, all $M_{n}$ encapsulate the same torus and all $M_{n+1/2}$ encapsulate the other. Also notice that any sphere which encapsulates both torus components of $\partial N$ is isotopic to $S$

\begin{lemma}\label{rotsym}
	Let $\Sigma\ll M_{k/2}$ . Then there exists a surface $\Sigma'\in\mathcal{I}(\Sigma)$ such that $\Sigma'\ll M_{k/2}$ and for which one of the components of $\Sigma'$ is rotationally symmetric.
\end{lemma}

\begin{proof}
	Let $\Sigma_1$ and $\Sigma_2$ be the disjoint spherical components of $\Sigma$ and let $D_1=\overline{\Sigma_1-M_{k/2}}$, $D_2=\overline{\Sigma_2-M_{k/2}}$, and $A=\overline{M_{k/2}-\Sigma}$. Let $B\subset N$ be the open ball with $\partial B=A\cup D_1\cup D_2$. 
	
	Because $\pi_1(M_{k/2})=0$ and $M_{k/2}$ is rotationally symmetric, we can find an isotopy $\psi$ of $\Sigma$ such that $\psi_1(M_{k/2})=M_{k/2}$ and $\psi_1(\Sigma)$ is another compression of $M_{k/2}$ such that $\psi(\partial D_1)$ and $\psi(\partial D_2)$ are rotationally symmetric circles in $M_n$. Without loss of generality, we will assume that $\partial D_1$ and $\partial D_2$ are rotationally symmetric. 

	Each component of $\Sigma$ encapsulates at least one torus component of $\partial N$. Otherwise, one of the components would bound a ball in $N$ and be homotopically trivial, contradicting the definition of compression Moreover, one component of $\Sigma$ encapsulates the other. If not, then each would encapsulate a separate torus, $B$ would lie in the region outside both spheres, and $\Sigma$ would have resulted from a compression of a surface which encapsulates both tori. Since $M_{k/2}$ encapsulates only one, this cannot be. Say $\Sigma_2$ encapsulates $\Sigma_1$. 

	If $\Sigma_2$ encapsulates only one torus, then the region between $\Sigma_2$ and $\Sigma_1$ is homeomorphically $S^2\times(-\epsilon,\epsilon)$. $B$ would have to lie inside this region. $M_{k/2}$ would then bound a ball in $N$ and be homotopically trivial, a contradiction. 
	
	Thus, $\Sigma_2$ encapsulates both tori and is isotopic to $S$. Therefore we can find an isotopy which fixes $\Sigma'-D_2$ and takes $\Sigma_2$ to a rotationally symmetric surface.
\end{proof}

\emph{Remark:} Since $D_1$ is homotopic to $A\cup D_2$, there is an isotopy of $\Sigma$ which takes $D_1$ to a rotationally symmetric disk in the euclidean $\epsilon$-neighborhood of $A\cup D_2$, thus we can, in fact, choose both components of $\Sigma'$ to be rotationally symmetric. The full result, however, is unnecessary in the following proof.

\begin{figure}
			\begin{center}
				\begin{overpic}[scale=1.5]{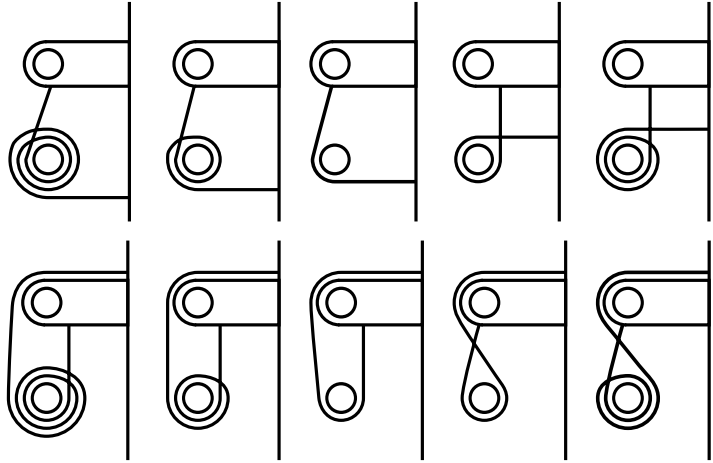}
				\end{overpic}
				\caption{Examples for $\gamma_D$}\label{slices}
			\end{center}
		\end{figure}

\begin{proof}[proof of proposition \ref{onlygammared}]
	Throughout this proof, if $\gamma$ is a curve, we will use $\gamma$ and $\im(\gamma)$ interchangeably.

	By lemma (\ref{rotsym}), there exists a compression $\Sigma'\in\mathcal{I}(\Sigma)$ for which one of the components is rotationally symmetric. Let $D_2'=\overline{\Sigma_2'-M_{k/2}}$ where $\Sigma_2'$ is the component of $\Sigma'$ which is rotationally symmetric. Now consider $\gamma_M=p(M_{k/2})$ and $\gamma_D=p(D_2')$ as subsets of $\Omega_{1/2}=im(p)$. Also let $L'=p(L)$. Since $M_{k/2}$ and $D_2'$ are embedded and rotationally symmetric, $\gamma_M$ and $\gamma_D$ are embedded 1-dimensional curves. Moreover, they are connected so we can consider $\gamma_M,\gamma_D:I\rightarrow \Omega_{1/2}$ as curves. 

		$\gamma_D$ has one endpoint in $L'$ and another in $\gamma_M$. The endpoints of $\gamma_M$ split $L'$ into three line segments. $\Sigma_2'$ is isotopic to $S$ so $\gamma_D$ must lie in the component of $\Omega_{1/2}-\gamma_M$ containing $p(S)$ and the endpoint of $\gamma_D$ in $L'$ must lie in one of the two outer subsegments of $L'-\gamma_M$.
	
	Let $\mathfrak{L}$ be the set of line segments in the component of $\Omega_{1/2}-\gamma_M$ containing $\gamma_D$ with one endpoint in $\gamma_M$ and the other in $L'$. Define the equivalence relation $\sim$ on $\mathfrak{L}$ as follows: $\delta\sim \delta'$ if there exists a homotopy of $\delta$ such that the endpoint of $\delta$ which is in $\gamma_M$ remains in $\gamma_M$ and the endpoint in $L'$ remains in $L'$ for the entire homotopy. The component of $\Omega_{1/2}-\gamma_M$ containing $\gamma_D$ is a homeomorphically a disk minus a sub-disk so $\pi_1$ of this component is $\Z$. Since the endpoint of $\gamma_D$ in $L'$ is contained in one of two path connected components, there is a one to two correspondence between $\Z$ and $\mathfrak{L}/\sim$. Representatives of these are shown in figure (\ref{slices}) with $M_0$ chosen for simplicity. We can see that there are only two possibilities for which $\gamma_D$ (and thus $D_2'$) is embedded, both of which correspond to the compression described.
\end{proof}

		\addtocounter{equation}{1}

To conclude this section, we will prove proposition (\ref{ordering}) which concerns compression and the infimum of areas in isotopy classes.

\begin{lemma}\label{orderlemma}
	Let $\Sigma'\ll\Sigma$ and $\epsilon>0$. There is an isotopy $\varphi$ of $\Sigma$ which fixes $\Sigma'$ and for which $|\varphi(\Sigma)|\leq|\Sigma'|+\epsilon$.
\end{lemma}
\begin{proof}
	Since $\overline{\Sigma'\Delta\Sigma}$ bounds a ball there is an isotopy which takes fixes $\Sigma'$ and which takes $\Sigma-\Sigma'$ to a region which approximates $\Sigma'-\Sigma$ with a thin tube connecting both disks. The tube can be made arbitrarily small, and the resulting surface has area less than $|\Sigma'|+\epsilon$
\end{proof}

\begin{proposition}\label{ordering}
	Let $\mu'=\inf_{\hat\Sigma\in\I(\Sigma')}|\hat\Sigma|$ and $\mu=\inf_{\hat\Sigma\in\I(\Sigma)}|\hat\Sigma|$. If $\Sigma'\ll\Sigma$, then $\mu'\geq\mu$.
\end{proposition}
\begin{proof}
	Let $\varphi$ be an isotopy of $\Sigma'$ such that $|\varphi(\Sigma')|<\mu'+\epsilon$. $\varphi(\Sigma)$ is an compression of $\varphi(\Sigma')$, so by lemma (\ref{orderlemma}), we can find $\hat\Sigma$ isotopic to $\varphi(\Sigma)$ with area less than $|\varphi(\Sigma')|+\epsilon$. So we have
	\begin{equation}\mu\leq|\hat\Sigma|\leq|\varphi(\Sigma')|+\epsilon\leq\mu'+2\epsilon\end{equation}

	Allowing $\epsilon$ to go to zero, we have our desired result.
\end{proof}

\section{Minimal Spheres}

Let $S_{k/2}$ be a MSY-minimizer for $M_{k/2}$. Let $S_{min}$ and $T_{min}$ be the MSY-minimizers for $S$ and $T$ respectively. We let $\mu_{k/2}=|S_{k/2}|$, $\sigma=|S_{min}|$, and $\tau=|T_{min}|$. In this section, we show that, under certain restrictions on the metric for $N$, we can guarantee a subsequence of $S_{k/2}$ which are isotopic to $M_{k/2}$.

\begin{lemma}\label{monotone}
	There exists an integer $m$ such that $S_{k/2}\in\I(M_{(k+m)/2}\cup |m|S_{min})$
\end{lemma}
\begin{proof}

	By proposition (\ref{Tseq}), we know that there is a finite sequence $M_{k/2}=\Sigma_0\gg\ldots\gg\Sigma_n$ such that $\Sigma_n\in\I(S_{k/2})$.  If $n=0$, then $m=0$ and we are done. Assume $n\geq1$. Proposition (\ref{onlygammared}) tells us that $\Sigma_k^1\in\I(M_{(k\pm 1)/2}\cup S)$ so if $n=1$, then $m=\pm 1$ and we are also done. 

	So assume $n>2$. First notice that, by proposition (\ref{ordering}), we have
	\begin{align}
		\mu_{k/2}	&	=	\inf_{\Sigma\in\I(\Sigma_0)}|\Sigma|
								\geq	\inf_{\Sigma\in\I(\Sigma_1)}|\Sigma|
								\geq\ldots
								\geq	\inf_{\Sigma\in\I(\Sigma_n)}|\Sigma|=\mu_{k/2}
	\end{align} 
	The last equality comes from the fact that $S_{k/2}$ is the MSY-minimizer for $M_{k/2}$. Therefore
	\begin{equation}\label{chainequality}
	\inf_{\Sigma\in\I(\Sigma_l)}|\Sigma|=\mu_{k/2} \text{ for all }0\leq l\leq n
	\end{equation}

	If $\Sigma_l\in\I(M_{i/2}\cup \hat S)$ where $l<n$, $i$ is a constant and $\hat S$ is some surface, then $\Sigma_{l+1}$ is isotopic to either $M_{(i\pm 1)/2}\cup S\cup \hat S$ or $M_{i/2}\cup \hat S'$ where $\hat S'\ll S$. Inductively we get that for $0\leq l\leq n$, $\Sigma_l\in\I(M_{l_i/2}\cup \hat S)$ for some $\hat S$ which is obtained from compressions of possibly multiple copies of $S$. Note that $l_{i+1}$ is equal to either $l_i\pm 1$ or $l_i$

	We claim that $l_i$ is either non-decreasing or non-increasing. To show this, assume that for some $0\leq i<j\leq n$ we have $l_i=l_j$. Assume there is an intermediate term $i<\alpha<j$ such that $l_\alpha\neq l_i$. without loss of generality we can assume that that $l_{i+1}\neq l_i$. Let $\Sigma_i=M\cup \tilde S$ where $M\in\I(M_{l_i/2})$, then we can say $\Sigma_{i+1}=\hat M\cup \hat S\cup \tilde S$ where $\hat M\cup \hat S\ll M$ and $\Sigma_j=M'\cup\hat S'\cup\tilde S'$ where $M'$, $\hat S'$, and $\tilde S'$ derive from a chain of compressions of $M$, $\hat S$ and $\tilde S$ respectively. To put it neatly
		\begin{equation}\Sigma_i=M\cup\tilde S\gg \hat M \cup \hat S\cup \tilde S\gg\ldots\gg M'\cup\hat S'\cup \tilde S'=\Sigma_j\end{equation}
	Using the fact that $M'\in\I(M_{l_i/2})$, we get
	\begin{align}
		\inf_{\Sigma\in\I(\Sigma_i)}|\Sigma|&=\mu_{l_i/2}+\inf_{\Sigma\in\I(\tilde S)}|\Sigma|	\\
		\inf_{\Sigma\in\I(\Sigma_j)}|\Sigma|&=\mu_{l_i/2}+\sigma+\inf_{\Sigma\in\I(\tilde S)}|\Sigma|
	\end{align}
	Since $\sigma\neq0$, this contradicts equation (\ref{chainequality}), so for all $i<\alpha<j$, $l_\alpha=l_j$ and $l_i$ is either non-increasing or non-decreasing. So we have $S_{k/2}\in\I(M_{(k+m)/2}\cup \hat S)$ for some $\hat S$. 
	
	To figure the exact nature of $\hat S$, consider again the fact that if $\Sigma_i\in\I(M_{l_i/2}\cup\tilde S)$ then $\Sigma_{i+1}$ is isotopic to either $M_{(l_i\pm1)/2}\cup S\cup \tilde S$ or to $M_{l_i/2}\cup\tilde S'$ where $\tilde S'\ll\tilde S$. Since every step of the chain is of one of these forms, we see that $\hat S$ must have resulted from compressions of the ``leftover'' components in the former compression (i.e., the components isotopic to $S$). Since this step happens exactly $|m|$ times, $\hat S$ must have resulted from a compression of a surface which was isotopic to $|m|S$, and since $S_{min}$ minimizes area in this isotopy class, we get that $\hat S\in\I(|m|S_{min})$.
\end{proof}

\begin{figure}
			\begin{center}
				\begin{overpic}[scale=1.5]{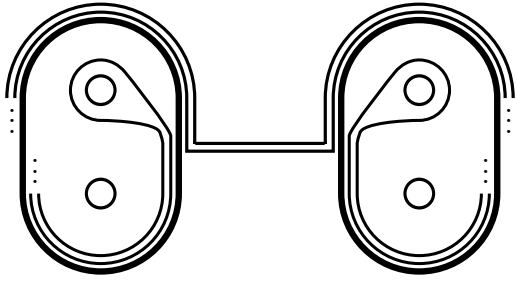}
					\put(.5,28){$\vdots$}
					\put(5,18){$\vdots$}
					\put(98,28){$\vdots$}
					\put(93.5,18){$\vdots$}
				\end{overpic}
				\caption{The cross section of a surface approximating $T_{min}$ (bold) and $M'$}\label{approximate}
			\end{center}
		\end{figure}

\begin{lemma}\label{bounds}
	There exists a constant $c$ such that $\mu_{k/2}\leq k\tau+c$
\end{lemma}
\begin{proof}
	$M_{k/2}$ is isotopic to a surface $M'$ which hooks around one of the torus components of $\partial N$, loops around in a region arbitrarily close in area to $kT$, and caps off with disks in the middle (see figure (\ref{approximate})). If we allow the hook and the disks to remain fixed, we find a surface with area less than $k\tau+c+\epsilon$ where $c$ is the area of the disk and the hook, and $\epsilon$ is arbitrarily small. Letting $\epsilon\rightarrow 0$, we get our desired result.
\end{proof}

		\addtocounter{equation}{1}

\begin{lemma}\label{isosequence}
	If $\tau<\sigma$, then there is a sequence $k_i$ such that $S_{k_i/2}\in\I(M_{k_i/2})$. Moreover $\mu_{k_i/2}$ is unbounded.
\end{lemma}
\begin{proof}
	By lemma (\ref{monotone}), for each $k$ there is an $m_k\in Z$ such that $S_{k/2}\in\I(M_{(k-m_k)/2}\cup |m_k|S)$. If there is no such sequence then $|k-m_k|<K$ for some large $K$ and all $k$, hence $\mu_{k/2}=\mu_{(k-m_k)/2}+m_k\sigma$. Then we have

\begin{align}
\begin{split}\label{adam}
	\sigma  &=    \lim_{k\rightarrow\infty}\frac{|k-K|\sigma}{k}
								\leq  \lim_{k\rightarrow\infty}\frac{\mu_{(k-m_k)/2}+m_k\sigma}{k}\\
					&=    \lim_{k\rightarrow\infty}\frac{\mu_k}{k}
								\leq  \lim_{k\rightarrow\infty}\frac{k\tau+c}{k}
								=     \tau
\end{split}
\end{align}

	Here the first inequality comes from the fact that $\{\mu_k:|k|<K\}$ is a finite set and the second inequality is lemma (\ref{bounds}). But this violates our assumption that $\tau<\sigma$, therefore such a sequence must exist.

	So assume there is a bound for area. By \cite{schoen1983}, there is a constant $C$ such that for small enough $r$ and $\rho\in(0,r]$, 
	\begin{equation}\sup_{B_{r-\rho}}|A_i|^2\leq\frac{C}{\rho^2}\end{equation}
	where $B_{r-\rho}$ is a ball of radius $r-\rho$ and $A_i$ is the second fundamental form $S_{k_i/2}$. $S_{k_i/2}$ is a sequence of stable, compact, connected, embedded minimal surfaces without boundary with a uniform bound on area and second fundamental form, therefore, by lemma (1) of \cite{dean2003} there is a compact connected embedded minimal surface without boundary, $\Sigma$, to which the sequence converges with finite multiplicity. By the maximum principle, we can see that for large $i$, $S_{k_i/2}$ are coverings of $\Sigma$, and since the $\mathbb{RP}^2$ does not embed into $\R^3$ (and thus $N$), $\Sigma$ must be a sphere and the multiplicity of convergence must be 1. Therefore, for large $i$, $\S_{k_i/2}\in\I(\Sigma)$, which contradicts the fact that $M_{k/2}$ are isotopically distinct.
\end{proof}
\begin{lemma}\label{metric}
	The set of metrics on $N$ for which $N$ is strictly mean convex and $\tau<\sigma$ is open and non empty
\end{lemma}
\begin{proof}
	Let $g$ be a metric on $N$ for which $\tau<\sigma$ and $\tilde g$ some other $C^2$ metric. Let $F:(N,g)\rightarrow(N,\tilde g)$ and let $k_F$ and $K_F$ be the minimum and maximum for $|df|$ respectively.  Let $\tilde \tau$ and $\tilde \sigma$ be the counterparts of $\tau$ and $\sigma$ in the new metric. Also, we let $|\Sigma|_g$ and $|\Sigma|_{\tilde g}$ be the area of $\Sigma$ in the respective metrics. Note that $k_f|\Sigma|_g\leq|\Sigma|_{\tilde g}\leq K_F|\Sigma|_g$. Then we have 
	\begin{align} \begin{split}
		\tilde \tau-\tilde\sigma	
			&=		\inf_{\Sigma\in\I(T)}|\Sigma|_{\tilde{g}} - \inf_{\Sigma\in\I(S)}|\Sigma|_{\tilde g}
			\leq 	K_F\inf_{\Sigma\in\I(T)}|\Sigma|_{g}-k_F \inf_{\Sigma\in\I(S)}|\Sigma|_{g}\\
			&=		K_F\tau-k_F\sigma	
	\end{split} \end{align}
	Choose a constant $\epsilon$ such that $x\tau<\sigma$ for all $x<1+\epsilon$. Then if $\frac{K_F}{k_F}<1+\epsilon$, the right-most part of the above equation is negative and we get $\tilde\tau<\tilde\sigma$. Thus the set of metrics for which $\tau<\sigma$ is open. Proposition (1) in \cite{dean2003} states that the set of metrics for which $N$ is mean convex is both open and non empty, so our set of metrics is indeed open. We need now only show the intersection of these two classes is non-empty.

	Let $x=(\rho,\theta,z)$ be cylindrical coordinates for $\R^3$ and for $r<1$, let
	\begin{align}\begin{split}
		T_r^+	&=	\{x\in\R^3:(\rho-1)^2+(z-\frac{r}{2})^2<\frac{r^2}{16}\}\\ 
		T_r^-	&=	\{x\in\R^3:(\rho-1)^2+(z+\frac{r}{2})^2<\frac{r^2}{16}\}\\
		B 		&=	\{x\in\R^3:\|x\|\leq 2\},\ N=B-(T_r^+\cup T_r^-)\\
		T_r		&=	\{x\in\R^3:(\rho-1)^2+z^2=r\}\\
		C			&=	\{x\in B:\rho<\frac{1}{2}\}
	\end{split}\end{align}
	So $N=B-(T_r^+\cup T_r^-)$ is as before. Also, note that $T_r\in\I(T)$ and $|T_r|=4\pi^2r$. Let $g$ be the euclidean metric on $N$. By proposition 1 in \cite{dean2003}, we can find $f\in C^2(N)$ which has support in an $\epsilon$ neighborhood of $\partial N$ for which $N$ with the metric $\tilde g=e^{2f}g$ is mean convex. Fix $\tilde g$ to be the metric for $N$.
	
	Notice that for any surface $\Sigma\in\I(S)$, $\Sigma\cap C$ contains at least two disks with area at least $\pi/4$. So, if $r<1/8\pi$, we have
	\begin{equation}\tau\leq 4\pi^2r<\pi/2\leq \sigma\end{equation}
	and $N$ with $\tilde g$ has the desired properties.
\end{proof}
\begin{proof}[Proof of Main Theorem]
	Let $M^3$ be a three-manifold. $N$ embeds into a three dimensional ball which can in turn be embedded into $M$, so we will consider $N$ to be a subset of $M$. The open set of metrics is that described in lemma (\ref{metric}), for each of which there is a sequence of stable, embedded two-spheres with arbitrarily large area as described in lemma (\ref{isosequence}).
\end{proof}

\bibliographystyle{alpha}
\bibliography{largesphere}

\end{document}